\documentclass[11pt,reqno]{article}
\usepackage{euscript,amstext,amsthm,amssymb}
\usepackage{amsmath}

\newcommand{\be}{\begin{equation}}
      \newcommand{\ee}{\end{equation}}
\newcommand{\ban}{\begin{eqnarray*}}
       \newcommand{\ean}{\end{eqnarray*}}
\newcommand{\ba}{\begin{eqnarray}}
       \newcommand{\ea}{\end{eqnarray}}

\font\BBb=msbm10 at 12pt
\renewcommand{\Bbb}[1]{\mbox{\BBb #1}}

\begin{document}
\newcommand{\pt}{\partial}
\newcommand{\lp}{\langle}
\newcommand{\rp}{\rangle}
\newcommand{\ra}{\rightarrow}
\newcommand{\rat}{\mapsto}
\newcommand{\lra}{\longrightarrow}
\newcommand{\ul}{\underline}
 \renewcommand{\o}[2]{\frac{#1}{#2}}
\newcommand{\hf}{\o{1}{2}}

\renewcommand{\qed}{\hspace*{\fill}\rule{3mm}{3mm}\quad}
 \newcommand{\Pf}{\noindent {\em Proof.} }
\newcommand{\Ex}{\noindent {\em Example.} }
\newcommand{\Rk}{\noindent {\em Remark.} }
\newcommand{\Def}{\noindent {\em Definition.} }
\newcommand{\heta}{\hat{\eta}}
\newcommand{\reta}{\overline{\eta}}
\newcommand{\ch}{\mbox{\rm ch}}
\newcommand{\sgn}{\mbox{\rm sgn}}
\newcommand{\tr}{\mbox{\rm tr}}
\newcommand{\Tr}{\mbox{\rm Tr}}
\newcommand{\DET}{\mbox{\rm DET}}
\newcommand{\hol}{\mbox{\rm hol}}
\newcommand{\End}{\mbox{\rm End}}
\newcommand{\spec}{\mbox{\rm spec}}
\renewcommand{\sf}{\mbox{\rm sf}}
\newcommand{\ind}{\mbox{\rm ind}}
\newcommand{\rk}{\mbox{\rm rk}}
\newcommand{\re}{\mbox{\rm Re}}
\newcommand{\Res}{\mbox{\rm Res}}
\newcommand{\HH}{{\rm H}}
\newcommand{\im}{\mbox{\rm Im}}
\newcommand{\inte}{\mbox{\rm int}}
\newcommand{\zn}{\stackrel{\circ}{\nabla}}
\newcommand{\zg}{\stackrel{\circ}{g}}

\renewcommand{\theequation}{\arabic{section}.\arabic{equation}}
\newcommand{\sect}[1]{\section{#1} \setcounter{equation}{0}}

\newtheorem{theo}{Theorem}[section]

\newtheorem{lem}[theo]{Lemma}
\newtheorem{prop}[theo]{Proposition}
\newtheorem{coro}[theo]{Corollary}

\title{\sc A Positive Mass Theorem for Spaces with Asymptotic SUSY Compactification}
\author{Xianzhe Dai}
\maketitle

\begin{abstract} We prove a positive mass theorem for spaces which asymptotically approach a flat Euclidean space times a
Calabi-Yau manifold (or any special honolomy manifold except the quaternionic K\"ahler). This is motivated by the very recent
work of Hertog-Horowitz-Maeda \cite{HHM}.
\end{abstract}

In general relativity, isolated gravitational systems are modelled
by asymptotically flat spacetimes. The spatial slices of such
spacetime are then asymptotically flat Riemannian manifolds. That
is, Riemannian manifolds $(M^n, g)$ such that $M=M_0 \cup
M_{\infty}$ with $M_0$ compact and $M_{\infty} \simeq {\Bbb R}^n -
B_R(0)$ for some $R>0$ so that in the induced Euclidean
coordinates the metric satisfies the asymptotic conditions \be
\label{afm} g_{ij}=\delta_{ij} + O(r^{-\tau}), \ \ \  \partial_k
g_{ij}=O(r^{-\tau -1}), \ \ \  \partial_k \partial_l
g_{ij}=O(r^{-\tau -2}). \ee Here $\tau >0$ is the asymptotic order
and $r$ is the Euclidean distance to a base point. The total mass
(the ADM mass) of the gravitational system can then be defined via
a flux integral \cite{ADM}, \cite{LP} \be \label{admm} m(g)=
\lim_{R \ra \infty} \frac{1}{4 \omega_n} \int_{S_R} (\partial_i
g_{ij} - \pt_j g_{ii} ) * dx_j . \ee Here $\omega_n$ denotes the
volume of the $n-1$ sphere and $S_R$ the Euclidean sphere with
radius $R$ centered at the base point.

If $\tau > \o{n-2}{2}$ and $n \geq 2$, then $m(g)$ is independent of the asymptotic coordinates $x_i$, and thus is an
invariant of the metric. The positive mass theorem \cite{SY1}, \cite{SY2}, \cite{SY3}, \cite{Wi} says that this total mass is
nonnegative
provided one has nonnegative local energy density.

\begin{theo}[Schoen-Yau, Witten] Suppose $(M^n, g)$ is an asymptotically flat  spin manifold of dimension $n \geq 3$ and of
order $\tau > \o{n-2}{2}$. If the scalar curvature $R \geq 0$, then $m(g) \geq 0$ and $m(g)=0$ if and only if $M={\Bbb R}^n$.
\end{theo}

\Rk The scalar curvature $R$ is the local energy density.

According to string theory \cite{CHSW}, our universe is really ten dimensional, modelled by  $M^{3,1} \times X$ where $X$ is
a Calabi-Yau 3-fold. This is the so called Calabi-Yau compactification, which motivates the spaces we now consider.

We consider the complete Riemannian manifolds $(M^n, g)$ such that $M=M_0 \cup M_{\infty}$ with $M_0$ compact and $M_{\infty}
\simeq ({\Bbb R}^k - B_R(0)) \times X$ for some $R>0$ and $X$ a compact simply connected Calabi-Yau manifold (or with any
other special honolomy except
$Sp(m)\cdot Sp(1)$) so that the metric on $M_{\infty}$ satisfies
\be \label{afmsu}
g=\zg + h, \ \ \  \zg=g_{{\mathbb R}^k} + g_X, \ \ \
h=O(r^{-\tau}),
\ \ \  \zn\! h =O(r^{-\tau -1}), \ \ \  \zn \zn \! h =O(r^{-\tau -2}).
\ee
Here $\zn$ is the Levi-Civita connection of $\zg$, $\tau >0$ is the asymptotical order. We will call $M$ a space with
asymptotic SUSY
compactification.

The mass for such a space is then defined by \be \label{admkkm}
m(g)= \lim_{R \ra \infty} \frac{1}{4 \omega_k vol(X)} \int_{S_R
\times X} ( \zn_{e_a^0}g_{ja} - \zn_{e_j^0} g_{aa} ) * dx_j
dvol(X). \ee Here $\{ e_a^0 \} = \{ \o{\pt}{\pt x_i}, f_{\alpha}
\}$ is an orthornormal basis of $\zg$, the $*$ operator is the one
on the Euclidean factor, the index $i, j$ run over the Euclidean
factor and the index $\alpha$ runs over $X$ while the index $a$
runs over the full index of the manifold. In fact, this reduces to
\[ m(g)= \lim_{R \ra \infty} \frac{1}{4 \omega_k vol(X)} \int_{S_R \times X} (\partial_i g_{ij} - \pt_j g_{aa} ) * dx_j
dvol(X). \]

\Rk If $\tau > \o{k-2}{2}$ and $k \geq 2$, then $m(g)$ is independent of the asymptotic coordinates.

Our main result is

\begin{theo}
Let $(M, g)$ be a complete spin manifold as above and the asymptotic order $\tau > \o{k-2}{2}$ and $k \geq 3$. If $M$ has
nonnegative scalar curvature, then $m(g) \geq 0$ and $m(g)=0$ if and only if $M={\Bbb R}^k \times X$.
\end{theo}

\Rk The result extends without change to the case with more than one end.  \\
\Rk Just like in the usual case, the restriction $k \geq 3$ has to do with getting the correct spin structure at the ends.
See section 5 for additional comments regarding the spin structures of the ends.

Our motivation comes from a very recent work of
Hertog-Horowitz-Maeda \cite{HHM} on the Calabi-Yau
compactifications. Using the existence result of Stolz \cite{S1},
\cite{S2} on metrics of positive scalar curvature, they
constructed classical configurations which has regions of
(arbitrarily large) negative energy density as seen from the four
dimensional perspective. This should be contrasted with the
positivity (nonnegativity) of the  total mass, as guaranteed by
Theorem 0.2. According to \cite{HHM}, physical consequences of the
negative energy density include possible violation of Cosmic
Censorship and new thermal instability.

The Lorentzian version of Theorem 0.2 will be discussed in a separate paper.
\newline

{\em Acknowledgement:} This work is motivated and inspired by the work of Gary Horowitz and his collaborators
\cite{HHM}. The author is indebted to Gary for sharing his ideas and for interesting discussions. The author would also like
to thank Is Singer for bringing them together and for useful discussion. Thanks are also due to Xiao Zhang for useful
comments.

\section{Manifolds with special holonomy}

For a complete Riemannian manifold $(M^n, g)$, the holonomy group $H\! ol(g)$ (with respect to a base point) is the subgroup
of
$O(n)$ generated by parallel translations along all loops at the base point. For simply connected irreducible  nonsymmetric
spaces, Berger has given a complete classification of possible holonomy groups, namely, $SO(n)$ which is the generic
situation, $U(m)$ (if $n=2m$) which is K\"ahler, $SU(m)$ for Calabi-Yau, $Sp(m)\cdot Sp(1)$ (if $n=4m$)
which is called quaternionic K\"ahler, $Sp(m)$ which is called hyper-K\"ahler,
$Spin(7)$ (if $n=8$), and $G_2$ (if $n=7$). Except the generic and K\"ahler cases, the rest are called special holonomy.

If a Riemannian manifold $(M, g)$ is spin, then one can consider spinors $\phi$ on $M$ which are sections of the spinor
bundle $S$. The Levi-Civita connection $\nabla$ of $g$ lifts to a connection of the spinor bundle, which will still be
denoted by the same notation. In fact, any metric connections lift in the same way. The Dirac operator
\[ D \phi = e_i \cdot \nabla_{e_i} \phi, \]
where $e_i$ is a local orthonormal basis of $M$ and $e_i \cdot$ is the Clifford multiplication. A spinor $\phi$ is parallel if
$\nabla \phi = 0$.

Implicitly, all these depend on the underlying spin structure, which is in one-to-one correspondence with elements of
$H^1(M, \ {\Bbb Z}_2)$ \cite{LM}. Thus, for simply connected manifolds, one has a unique spin structure. It seems that the
issue of spin structure in this context is a subtle one, deserving further study. (See also section 5.)

All manifolds with special holonomy, with the exception of the quaternionic K\"ahler ones, carry nonzero parallel spinor. In
fact, one has the
following theorem of McKenzie Wang \cite{Wa}.

\begin{theo} \label{wang}
Let $(M, g)$ be a complete, simply connected, irreducible Riemannian spin manifold and $N$ be the dimension of parallel
spinors. Then $N > 0$ if and only if the holonomy group is one of $SU(m)$, $Sp(m)$, $Spin(7)$, $G_2$.
\end{theo}
\Rk Wang \cite{Wa} actually characterizes each special holonomy by the number of parallel spinors.

\Rk Manifolds with parallel spinors are called supersymmetric (SUSY) in physics literature.

\section{Proof of Theorem 0.2}

Our proof is an extension of Witten's spinor proof \cite{Wi}. Here we follow the idea of Anderson and Dahl \cite{AnD} and use
the following alternative formula for the Lichnerowicz formula.

\begin{lem} Given a spinor $\phi$ on a Riemannian spin manifold, define a $1$-form $\alpha$ via
\[ \alpha (X)= \lp (\nabla_X + X \cdot D)\phi, \ \phi \rp. \]
Then
\[ div\,\alpha = \o{R}{4} |\phi|^2 + |\nabla \phi|^2 - |D \phi|^2 .\]
\end{lem}

\Pf Choose an orthonormal basis $e_a$ such that $\nabla e_a = 0$ at the given point. Then (Einstein summation enforced)
\ba div\,\alpha & = & (\nabla_{e_a} \alpha)(e_a) = e_a (\alpha(e_a)) \nonumber \\
 & = & \lp (\nabla_{e_a} + e_a \cdot D)\phi, \ \nabla_{e_a}\phi \rp +  \lp \nabla_{e_a}(\nabla_{e_a} + e_a \cdot D)\phi, \
\phi \rp \nonumber \\
 & = & |\nabla \phi|^2 - |D \phi|^2 +  \lp (\delta_{ab} + e_a \cdot e_b \cdot) \nabla_{e_a} \nabla_{e_b} \phi, \ \phi \rp.
\nonumber
\ea
The last term is just
\[ \lp \hf [e_a \cdot, \  e_b \cdot] \nabla_{e_a} \nabla_{e_b} \phi, \ \phi \rp = \lp \o{1}{4} [e_a \cdot, \  e_b \cdot]
R(e_a,
e_b) \phi, \ \phi \rp = \o{R}{4} |\phi|^2 \]
by the usual calculation as in the Lichnerowicz formula \cite{LM}.  \qed

Therefore, for any compact domain $\Omega \subset M$,
\ba \label{beq}
\int_{\Omega} [\o{R}{4} |\phi|^2 + |\nabla \phi|^2 - |D \phi|^2 ]\, dvol(g) & = & \int_{\pt \Omega} \sum \lp (\nabla_{e_a} +
e_a \cdot D)\phi, \ \phi \rp \, \inte(e_a) \, dvol(g), \nonumber \\
& = & \int_{\pt \Omega} \sum \lp (\nabla_{\nu} + \nu \cdot D)\phi, \ \phi \rp \, dvol(g|_{\pt \Omega})
\ea
where $e_a$ is an orthonormal basis of $g$ and $\nu$ is the unit outer normal of $\pt \Omega$. Also, here $\inte(e_a)$ is the
interior multiplication by $e_a$.

In particular, for a harmonic spinor $\phi$, i.e., $D\phi =0$, the
left hand side of (\ref{beq}) will be nonnegative provided $R \geq
0$. On the other hand, if the harmonic spinor $\phi$ can be chosen
so that it is asymptotic to a parallel spinor at infinity and we
choose the domain $\Omega$ so that $\pt \Omega = S_R \times X$,
then we will show that the right hand side of (\ref{beq})
converges to the mass (up to a positive normalizing constant).
Thus, for the first part of our theorem, we are left with two
tasks. First, we need to show the existence of harmonic spinors
which are asymptotic to a parallel spinor. Second, we need to show
that the limit of the boundary term converges to the mass. The
existence of the harmonic spinor is dealt with in section 4 (Lemma
4.1) after the necessary analysis in the next section and the
computation of the limit of the boundary term is also left to
Section 4 (Lemma 4.2).

We now continue with the proof of the rigidity. If $m(g)=0$, then it follows that $\phi$ is a (nonzero) parallel spinor on
$M$. This implies that $M$ is Ricci flat, as
\[ e_a \cdot R(e_a, X) \phi = - \hf Ric(X) \, \phi .\]
Thus, we are in a position to use the splitting theorem of Cheeger-Gromoll \cite{CG}. To find lines in $M$, we start with
sequences of pairs of points $p_i, q_i$ in $ M_{\infty} \simeq ({\Bbb R}^k - B_R(0)) \times X$. When $R$ is sufficiently
large, one can choose $p_i, q_i$ so that their distance is comparable to their Euclidean distance. It follows that one can
construct a line in $M$ this way. Similarly, we can construct $k$ lines in $M$ that are almost perpendicular to each other.
It follows that $M= {\Bbb R}^k \times X$. \qed

\section{Fibered boundary calculus}

We will use the fibered boundary calculus of Melrose-Mazzeo \cite{MM} (and further developed by Boris Vaillant in his thesis
\cite{V} and in \cite{HHMa}) to solve for the harmonic spinor with the correct asymptotic behavior.

The change of variable $r=\o{1}{x}$ makes metric into what is called fibered boundary metric, which is defined in the more
general setting as follows.

Consider a complete noncompact Riemannian manifold $(M, g)$. Assume that $M$ has a compactification $\bar{M}$ such that $\pt
\bar{M}$ comes with a fibration structure $F \ra \pt \bar{M} \stackrel{\pi}{\longrightarrow} B$. Moreover, in a neighborhood
of the boundary $\pt \bar{M}$, the metric $g$ has the form
\be \label{fbm}
g = \o{dx^2}{x^4} + \o{\pi^*(g_B)}{x^2} + g_F
\ee
where $x$ is a defining function of the boundary, i.e., $x=0$ on $\pt \bar{M}$ and $dx \not= 0$ on the boundary. Also, $g_B$
is a metric on the base $B$, $g_F$ is a family of fiberwise metrics.

Thus, in the setting of spaces with asymptotic SUSY compactification, one has a trivial fibration $S^{k-1} \times X$ and
$x=\o{1}{r}$.

We will use the notation $M$, $\bar{M}$, and $\pt M$, $\pt \bar{M}$ interchangeably. For a manifold with boundary, the Lie
algebra of $b$-vector fields consists of vector fields tangent to the boundary
\[ {\cal V}_b(M)= \{ V \ | \ V \ \mbox{is tangent to the boundary} \ \pt M \}  \]
The Lie algebra of vector fields associated with the fibered boundary metric is
\be \label{fbvf}
{\cal V}_{fb} = \{ V \in {\cal V}_b(M)\ | \ V \ \mbox{is tangent to the fibers $F$ at} \ \pt M, \ Vx=O(x^2) \} .
\ee
If $y$ is local coordinates of $B$ and $z$ is local coordinates of $F$, then ${\cal V}_{fb}$ is spanned by $x^2 \pt_x, \
x\pt_y, \pt_z$. The fibered boundary vector fields ${\cal V}_{fb}$ generate the ring of fibered boundary differential
operators. The Dirac operator $D$ associated to the fibered boundary metric is such a fibered boundary differential operator
of first order.

Define the $L^2$ and Sobolev spaces as follows.
\[ L^2(M, S) = L^2 (M, S; dvol(g)) = L^2(M, S, \o{dxdydz}{x^{2+l}}) \]
if $\dim B=l$.
\[ L^{p, 2}(M, S) = \{ \ \phi \in L^2(M, S) \ | \ \nabla_{V_1} \cdots \nabla_{
V_j} \phi \in L^2(M, S), \ \forall j \leq p, \ V_i \in {\cal V}_{b} \ \}. \]

For $\gamma \in {\Bbb R}$, the space of conormal sections of order $\gamma$ is defined to be
\[ {\cal A}^{\gamma}(M, S) = \{ \ \phi \in C^{\infty}(M, S) \ | \ \nabla_{V_1} \cdots \nabla_{V_j} \phi | \leq C x^{\gamma},
\ \forall j, \ V_i \in {\cal V}_{b} \ \}, \]
while the space of polyhomogeneous sections is
\[ {\cal A}^{*}_{phg}(M, S) = \{ \ \phi \in {\cal A}^*(M, S) \ | \  \phi \sim \sum_{\re \gamma_j \ra \infty} \sum_{k=0}^{N_j}
\psi_{jk} x^{\gamma_j} (\log x)^k,  \  \psi_{jk} \in C^{\infty}(\pt M, S) \ \}. \]
Here the expansion is the usual asymptotic expansion, uniform with all the derivatives.
We usually specify all possible pair $(\gamma_j, \, N_j)$ that can appear in the expansion and the collection of $(\gamma_j,
\, N_j)$ is called the index set.

Assume that $\ker D_F$ has constant dimension so it forms a vector bundle on the base $B$. Let $\Pi_0$ be the orthogonal
projection onto $\ker D_F$ and $\Pi_{\perp}= I –- \Pi_0$. The following is a summary of the
results developed in \cite{MM}, \cite{V}, \cite{HHMa}.

\begin{theo} \label{fbc}
Suppose that $a$ is not an indicial root of $\Pi_0 x^{-1} D \Pi_0$. Then
\[ D:\, x^a L^{1,2}(M, S) \ra x^{a+1} \Pi_0 L^{2}(M, S) \oplus x^a \Pi_{\perp} L^{2}(M, S) \]
is Fredholm. If $D\phi=0$ for $\phi \in x^a L^{2}(M, S)$, then $\phi$ is polyhomogeneous with exponents in its expansion
determined by the indicial roots of
$\Pi_0 x^{-1}D \Pi_0$ and truncated at $a$. If $D\xi = \psi$ for $\psi \in {\cal A}^a(M, S)$ and $\xi \in x^{c-1} \Pi_0
L^{1,2}(M, S) \oplus x^c
\Pi_{\perp} L^{1, 2}(M, S)$ and $c < a$, then $\xi \in \Pi_0 {\cal A}^I_{phg}(M, S) + {\cal A}^a(M, S)$.
\end{theo}

For the precise definition of the indicial root, and in particular, the indicial
root of $\Pi_0 x^{-1} D \Pi_0$, we refer the reader to \cite{MM}, \cite{HHMa}. For our purpose, we only note that it is a
discrete set.

\Rk Strictly speaking, only $\zg$ is a fibered boundary metric in the pure sense but it is easy to see that the result
generalize to the metric $g$. In any case, the metric perturbation produces only a lower order term (Cf. section 4).

\begin{lem} \label{iso}
If $R\geq 0$ and $a> \o{k-2}{2}$ is not an indicial root, then
\[ D:\, x^a L^{1,2}(M, S) \ra x^{a+1} \Pi_0 L^{2}(M, S) \oplus x^a \Pi_{\perp} L^{2}(M, S) \]
is an isomorphism.
\end{lem}
\Pf We first see that it is injective. If $D\phi=0$ for $\phi \in x^a L^{2}(M, S)$, then by Theorem \ref{fbc}, $\phi \in
{\cal A}^a_{phg}(M, S)$. Now, from (\ref{beq}),
\[ \int_{\Omega} [ |\nabla \phi|^2 + \o{R}{4} |\phi|^2 ] dvol = \int_{\pt \Omega} \lp \nabla_{\nu}\phi, \, \phi
\rp \, dvol(\pt \Omega). \]
By taking $\Omega$ so that $\pt \Omega = S_r \times X$ and $r \ra \infty$ we see that the right hand side goes to zero since
$\phi \in {\cal A}^a_{phg}(M, S)$ and $a> \o{k-2}{2}$. It follows then by the assumption $R \geq 0$ that $\phi$ is parallel
and
hence zero.

Now, if $\omega$ is in the cokernel of $D$, then, by the Fredholm property, $\omega \in x^{a+1} \Pi_0 L^{2}(M, S) \oplus x^a
\Pi_{\perp} L^{2}(M, S)$ and $\omega$ is a weak solution of Dirac equation:
\[ \lp \omega, \ D \xi \rp =0, \ \forall \xi \in x^a L^{1,2}(M, S). \]
It follows by the regularity part of Theorem \ref{fbc}, $\omega \in {\cal A}^a_{phg}(M, S)$. Therefore the same argument as
above shows $\omega =0$.    \qed

\section{Computation of the mass}

Recall that $g=\zg + h$ with $\zg=g_{{\mathbb R}^k} + g_X$ and $h=O(r^{-\tau})$,
$\zn\! h =O(r^{-\tau -1})$, $\zn \zn\! h =O(r^{-\tau -2})$. Let $e^0_a$ be the orthonormal basis of $\zg$ which consists of
$\o{\pt}{\pt x_i}$ followed by an orthonormal basis $f_{\alpha}$ of $g_X$. Orthonormalizing $e^0_a$ with respect to $g$ gives
rise an orthonormal basis $e_a$ of $g$. Moreover,
\be \label{onbasy}
e_a=e^0_a - \hf h_{ab} e^0_b + O(r^{-2\tau}).
\ee
This gives rise to a gauge transformation
\[ A: \ SO(\zg) \ni e^0_a \ra e_a \in SO(g) \]
which identifies the corresponding spin groups and spinor bundles.

To compare $\nabla$ and $\zn$, in particular their lifts to the spinor bundles, one introduces a new connection $\nabla^0=A
\circ \zn \circ A^{-1}$. This connection is compatible with the metric $g$ but has a torsion
\be \label{tor}
T(X, Y)=\nabla^0_X Y - \nabla^0_Y X – [X, Y] = -(\zn_X A) A^{-1} Y + (\zn_Y A) A^{-1} X .
\ee
The difference of $\nabla$ and $\nabla^0$ is then expressible in terms of the torsion
\be \label{diffcs}
2\lp \nabla^0_X Y - \nabla_X Y, Z \rp = \lp T(X, Y), Z\rp - \lp T(X, Z), Y \rp - \lp T(Y, Z), X \rp,
\ee
where we use the metric $g$ for the inner product $\lp \ , \ \rp$.

Since $\nabla$ and $\nabla^0$ are both $g$-compatible, their induced connections on the spinor bundle differ by
\be \label{diffsc}
\nabla_{e_a} - \nabla^0_{e_a} = -\o{1}{4} \sum_{b, c} (\omega_{bc}(e_a) - \stackrel{\circ}{\omega}_{bc}(e_a))e_b e_c,
\ee
where $e_b, e_c$ act on the spinors by the Clifford multiplication and the connection $1$-forms
\[ \omega_{bc}(e_a)=\lp \nabla_{e_a} e_b, e_c \rp, \ \ \ \stackrel{\circ}{\omega}_{bc}(e_a)=\lp \zn_{e_a}\! e_b, e_c \rp. \]
From (\ref{diffcs}) and (\ref{tor}) we obtain
\be \label{fundeq}
\nabla_{e_a} - \nabla^0_{e_a} = \o{1}{8} \sum_{b\not= c} (\zn_{e_b} g_{ac} - \zn_{e_c} g_{ab}) e_b e_c + O(r^{- 2\tau -1})
\ee
for the difference of the two connections acting on spinors.

\begin{lem} \label{hspinor}
There exists a harmonic spinor on $(M, \ g)$ which is asymptotic
to a parallel spinor at infinity.
\end{lem}
\Pf Our manifold $M= M_0 \cup M_{\infty}$ with $M_0$ compact and $M_{\infty} \simeq ({\Bbb R}^k - B_R(0)) \times X$. Since $k
\geq 3$ and $X$ is simply connected, the end $M_{\infty}$ is also simply connected, and therefore has a unique spin structure
coming from the product of the restriction of the spin structure on ${\Bbb R}^k$ and the spin structure on $X$.

Now pick a unit norm parallel spinor $\psi_0$ of $({\Bbb R}^k, g_{{\mathbb R}^k})$ and a unit norm parallel spinor $\psi_1$
of $(X, g_X)$. Then $\phi_0=A(\psi_0 \otimes \psi_1)$ defines a spinor of $M_{\infty}$.  We extend $\phi_0$ smoothly inside.
Then
$\nabla^0 \phi_0 =0$ outside the compact set. Thus, it follows from (\ref{fundeq}) that
\be \label{apsp}
\nabla \phi_0 = O(r^{-\tau -1}).
\ee

We now construct our harmonic spinor by setting $\phi=\phi_0 + \xi$ and solve $D\xi=-D\phi_0 \in O(r^{-\tau -1})$. By using
Lemma \ref{iso}, adjusting $\tau$ slightly if necessary so that it is not one of the indicial root, we have a solution $\xi
\in O(r^{-\tau})$.  \qed

\begin{lem} \label{comt}
For the harmonic spinor $\phi$ constructed above, we have
\[ \lim_{R \ra \infty} \int_{S_R \times X} \sum \lp (\nabla_{e_a} + e_a \cdot D)\phi, \ \phi \rp \, \inte(e_a) \, dvol(g) =
\omega_k vol(X) m(g). \]
\end{lem}

\Pf By (\ref{beq}),
\[ \int_{S_R \times X} \sum \lp (\nabla_{e_a} + e_a \cdot D)\phi, \ \phi \rp \, \inte(e_a) \, dvol(g)  = \re \int_{S_R \times
X} \sum \lp (\nabla_{e_a} + e_a \cdot D)\phi, \ \phi \rp \, \inte(e_a) \, dvol(g) . \]
Now,
\ba \lp (\nabla_{e_a} + e_a \cdot D)\phi, \ \phi \rp & = & \lp \hf [ e_a \cdot, e_b \cdot ]\, \nabla_{e_b}\phi, \ \phi \rp
\nonumber \\
& = & \lp \hf [ e_a \cdot, e_b \cdot ]\, \nabla_{e_b}\phi_0, \ \phi_0 \rp + \lp \hf [ e_a \cdot, e_b \cdot ]\,
\nabla_{e_b}\phi_0, \ \xi \rp  \nonumber  \\
& & + \ \lp \hf [ e_a \cdot, e_b \cdot ]\, \nabla_{e_b}\xi, \ \phi_0 \rp + \lp \hf [ e_a \cdot, e_b
\cdot ]\, \nabla_{e_b}\xi, \ \xi \rp.
\ea
The second term and the last term are $O(r^{- 2\tau -1})$ and therefore contribute nothing in the limit. For the third term,
one notice that if $\beta$ is the $n-2$ form
\[ \beta= \lp [ e_a \cdot, e_b \cdot ]\, \phi, \ \psi \rp \, \inte(e_a) \, \inte(e_b) \, dvol(g) \]
(Einsterin summation here and below), then
\ba d \, \beta & = & -2 \left( \lp [ e_a \cdot, e_b \cdot ]\, \nabla_{e_b} \phi, \ \psi \rp \, \inte(e_b) \, dvol(g) +
\lp [ e_a \cdot, e_b \cdot ]\, \phi, \ \nabla_{e_b} \psi \rp \, \inte(e_b) \, dvol(g) \right) \nonumber \\
& = & - 4 \left( \lp [ e_a \cdot, e_b \cdot ]\, \nabla_{e_b} \phi, \ \psi \rp \, \inte(e_b) \, dvol(g) -
\lp \phi, \ [ e_a \cdot, e_b \cdot ]\, \nabla_{e_b} \psi \rp \, \inte(e_b) \, dvol(g) \right) \nonumber \\
\ea
which yields
\[ \int_{\pt \Omega} \lp [ e_a \cdot, e_b \cdot ]\, \nabla_{e_b} \phi, \ \psi \rp \, \inte(e_b) \, dvol(g) = \int_{\pt \Omega}
\lp \phi, \ [ e_a \cdot, e_b \cdot ]\, \nabla_{e_b} \psi \rp \, \inte(e_b) \, dvol(g). \]
It follows then that the third term is similarly dealt with as the second. Thus the only contribution is coming from the
first term, for which we note that
\ban \lefteqn{ \lp \hf [ e_a \cdot, e_b \cdot ]\, \nabla_{e_b}\phi_0, \ \phi_0 \rp } \\
& = & \lp \hf [ e_a \cdot, e_b \cdot ] (\nabla_{e_b} - \nabla^0_{e_b}) \phi_0, \ \phi_0 \rp \\
& = & \o{1}{16} \sum_{c\not= d} (\zn_{e_c} g_{bd} - \zn_{e_d} g_{bc}) \lp
[e_a \cdot, \, e_b \cdot ] \, e_c \cdot e_d \cdot \phi_0, \ \phi_0 \rp + O(r^{- 2\tau -1 })
\ean
by (\ref{fundeq}). Now
\ban \lefteqn{ \o{1}{16} \sum_{c\not= d} (\zn_{e_c} g_{bd} - \zn_{e_d} g_{bc}) \lp
[e_a \cdot, \, e_b \cdot ] \, e_c \cdot e_d \cdot \phi_0, \ \phi_0 \rp } \\
& = & \o{1}{8} \sum_{c\not= d} (\zn_{e_c} g_{bd} - \zn_{e_d} g_{bc}) \lp
e_a \cdot e_b \cdot  e_c\cdot  e_d \cdot \phi_0, \ \phi_0 \rp  \\
&  &  + \ \o{1}{8} \sum_{c\not= d} (\zn_{e_c} g_{ad} - \zn_{e_d} g_{ac}) \lp
e_c \cdot, \, e_d \cdot \phi_0, \ \phi_0 \rp  \\
& = & \o{1}{8} \sum_{c\not= d} \zn_{e_c} g_{bd} \lp
e_a \cdot e_b \cdot  e_c\cdot  e_d \cdot \phi_0, \ \phi_0 \rp + \o{1}{8} \sum_{c\not= d} \zn_{e_d} g_{bb} \lp
e_a \cdot e_d \cdot \phi_0, \ \phi_0 \rp  \\
&  & + \ \o{1}{8} \sum_{c\not= d} (\zn_{e_c} g_{bd} - \zn_{e_d} g_{bc}) \lp  e_c\cdot  e_d \cdot \phi_0, \ \phi_0 \rp \\
& = & \o{1}{8} \sum_{c\not= d} \zn_{e_c} g_{bb} \lp
e_a \cdot e_c \cdot \phi_0, \ \phi_0 \rp  +– \o{1}{4} \sum_{c\not= d} \zn_{e_b} g_{bd} \lp e_a \cdot e_d \cdot \phi_0, \
\phi_0 \rp \\
&  &  +  \ \o{1}{8} \sum_{c\not= d} \zn_{e_d} g_{bb} \lp
e_a \cdot e_d \cdot \phi_0, \ \phi_0 \rp + \o{1}{8} \sum_{c\not= d} (\zn_{e_c} g_{bd} - \zn_{e_d} g_{bc}) \lp e_c\cdot  e_d
\cdot \phi_0, \ \phi_0 \rp\\
\ean For the last equality, we use $ e_c\cdot  e_d \cdot = \hf [
e_c\cdot, \,  e_d \cdot]$ for $c\not= d$, and $[ e_c\cdot, \,  e_d
\cdot]$ skew-hermitian to see that its real part is zero. Finally,
one uses $ e_a\cdot  e_d \cdot = \hf [ e_a\cdot, \,  e_d \cdot] -
\delta_{ad}$ and the skew-hermitian property of the commutators to
obtain
\ban \lefteqn{ \re \left( \lp \hf [ e_a \cdot, e_b \cdot ]\, \nabla_{e_b}\phi_0, \ \phi_0 \rp \right) } \\
 & = & \o{1}{4} ( \zn_{e_b}g_{ab} - \zn_{e_a} g_{bb} ) |\phi_0|^2 + O(r^{-2\tau -1}). \nonumber
\ean

This yields
\ban \lefteqn{ \lim_{R \ra \infty} \int_{S_R \times X} \sum \lp (\nabla_{e_a} + e_a \cdot D)\phi, \ \phi \rp \, \inte(e_a) \,
dvol(g) } \\
& = & \ \lim_{R \ra \infty} \int_{S_R \times X} \o{1}{4} ( \zn_{e_b}g_{ab} - \zn_{e_a} g_{bb} ) |\phi_0|^2 \inte(e_a) \,
dvol(g).
\ean
To see that this reduces to the definition of the mass, we first note that one can replace $e_a$ by $e_a^0$ in the integrand
on the right hand side, producing
only an error of $O(r^{-2\tau -1})$, then replace $dvol(g)$ by $dx dvol_X$ with a similar error term. \qed

\section{Negative energy solutions in Kaluza-Klein theory}

It was observed by Witten that positive energy theorems do not extend immediately to Kaluza-Klein theory \cite{Wi2}. He
observed that there are two zero energy solutions on a space asymptotic to $M_4 \times S^1$ which should lead to
perturbatively negative energy solutions. The explicit negative energy solutions were constructed later in \cite{BP},
\cite{BH}. The following example is from \cite{BH}.

The analytically continued Reissner-Nordstr\"om metric
\[ ds^2= (1-\o{2m}{r}-\o{q^2}{r^2}) d\theta^2 + (1-\o{2m}{r}-\o{q^2}{r^2})^{-1} dr^2 + r^2 d\Omega^2, \]
where $r \geq r_+=m + \sqrt{m^2+q^2}$, $\theta \in {\Bbb R}/\! \o{2\pi r_+^2}{r_+ -m} {\Bbb Z}$ and $d\Omega^2$ is the
standard metric on the $2$-sphere. This is a scalar flat metric on ${\Bbb R}^2 \times S^2$ and asymptotic to ${\Bbb R}^3
\times S^1$ at infinity. The mass can be computed via (\ref{admkkm}),  which is
\be \label{nmkk}
m(g)= \hf m \o{r_+ -m}{2\pi r_+^2}.
\ee
For fixed asymptotic geometry, i.e., fixed circle size $\o{2\pi r_+^2}{r_+ -m}= l$, this can be made arbitrarily negative if
one takes $m <0$ sufficiently large, while $q\not= 0$ is chosen appropriately (which will necessarily be large as well).

The reason here is that the end ${\Bbb R}^3 \times S^1$, and in particular, $S^1$ has the wrong spin structure! Recall that
$S^1$ has two spin structures
which correspond to the trivial double cover of $S^1$ and the nontrivial double cover of $S^1$. Here, since $S^1$ bounds the
disk inside, it has the spin structure corresponding to the nontrivial double cover. It therefore has no parallel spinor.

\end{document}